\documentclass[12pt, a4paper]{amsart}

\usepackage{amsthm,amsmath,amssymb}
\usepackage{graphicx}
\usepackage{tikz}
\usepackage{url}
\usepackage{a4wide}
\usepackage{microtype}
\usepackage{lmodern}
\usepackage[T1]{fontenc}

\usepackage[pdftex, plainpages = false, pdfstartview=FitH,
    pdfpagelabels,  pdfpagelayout = SinglePage, 
    bookmarks = true, bookmarksopen = true, 
    bookmarksnumbered = true, linktocpage,
    colorlinks = true, 
    linkcolor = blue, urlcolor  = blue, 
    citecolor = red,  anchorcolor = green, 
    hyperindex = true, hyperfigures]{hyperref} 

\tikzset{/tikz/ampersand replacement=\col}
\usetikzlibrary{calc,angles,positioning,intersections,quotes,decorations.markings}

\theoremstyle{plain}
\newtheorem{thm}{Theorem}
\newtheorem{lemma}[thm]{Lemma}
\newtheorem{corollary}[thm]{Corollary}
\newtheorem{proposition}[thm]{Proposition}

\theoremstyle{definition}
\newtheorem{defn}{Definition}
\newtheorem{example}[defn]{Example}
\newtheorem{conjecture}[defn]{Conjecture}

\theoremstyle{remark}
\newtheorem{remark}[thm]{Remark}

\newcommand{\bt}[1]{\begin{thm}\label{#1}}
\newcommand{\bc}[1]{\begin{corollary}\label{#1}}
\newcommand{\bl}[1]{\begin{lemma}\label{#1}}
\newcommand{\bp}[1]{\begin{proposition}\label{#1}}
\newcommand{\be}[1]{\begin{example}\label{#1}}
\newcommand{\bd}[1]{\begin{defn}\label{#1}}
\newcommand{\br}[1]{\begin{remark}\label{#1}}
\newcommand{\bx}[1]{\begin{exercise}\label{#1}}
\newcommand{\bcon}[1]{\begin{conjecture}\label{#1}}

\newcommand{\et}{\end{thm}}
\newcommand{\ec}{\end{corollary}}
\newcommand{\el}{\end{lemma}}
\newcommand{\ep}{\end{proposition}}
\newcommand{\ee}{\end{example}}
\newcommand{\ed}{\end{defn}}
\newcommand{\exc}{\end{exercise}}
\newcommand{\er}{\end{remark}}
\newcommand{\econ}{\end{conjecture}}

\newcommand{\bpr}{\begin{proof}}
\newcommand{\epr}{\end{proof}}
\newcommand{\rank}{\operatorname{rank}}
\newcommand{\ol}{\overline}

\newcommand{\tld}{\widetilde}
\newcommand{\lbu}{L^{\bullet}}
\newcommand{\ma}{M_{\mathcal{A}}}

\newcommand{\ua}{U_{\mathcal{A}}}
\newcommand{\tdhat}[1]{\widehat{T^{#1}}}
\newcommand{\redh}{\widetilde{H}}

\def\A  {\mathcal{A}}

\def \D {\mathcal{D}}

\def\R{\mathbb{R}}

\def\C{\mathbb{C}}
\def\Z{\mathbb{Z}}

\newcommand{\ctor}[1]{T^{#1}_{\C}}

\def\ds{\displaystyle}

\title{The Goresky-MacPherson formula for toric arrangements}
\author[P\, Deshpande]{Priyavrat Deshpande}
\thanks{The author is partially funded by a grant from Infosys Foundation}
\address{Chennai Mathematical Institute, India}
\email{pdeshpande@cmi.ac.in}
\date{}

\begin{document}

\begin{abstract}
A subspace arrangement is a finite collection of affine subspaces in $\R^n$.  
One of the main problems associated to arrangements asks up to what extent the topological invariants of the union of these spaces, and of their complement are determined by the combinatorics of their intersection.
The most important result in this direction is due to Goresky and MacPherson.
As an application of their stratified Morse theory they showed that the additive structure of the cohomology of the complement is determined by the underlying combinatorics.

In this paper we consider toric arrangements; a finite collection of subtori in $(\C^{\ast})^l$. 
The aim of this paper is to prove an analogue of the Goresky-MacPherson's theorem in this context. 
When all the subtori in the arrangement are of codimension-$1$ we give an alternate proof of a theorem due to De Concini and Procesi.
\end{abstract}
\maketitle


\section{Introduction}\label{intro}
A subspace arrangement in $\R^n$ is a finite collection 
$\A = \{A_1,\dots, A_k\} $
of affine subspaces such that if $A_i\subset A_j$ then $i=j$.
Two important spaces associated with an arrangement are the \emph{link} $V_{\A}$ which is the union of all the subspaces and the \emph{complement} $M_{\A}$ which is the complement of the link.
The formal data associated with an arrangement is the poset $L_{\A}$ of all non-empty intersections of members of $\A$ and the dimensions of these intersections. 
The poset $L_{\A}$ is ordered by reverse inclusion and the  ambient space $\hat{0}:= \R^n$ is its least element.
An arrangement of hyperplanes, now a subject in itself, is an important class of subspace arrangements.

A prominent research direction, in this area, is to try and express topological invariants of the link and those of the complement in terms of the combinatorics of the intersections. 
To this effect a seminal result is due to Goresky and MacPherson \cite[III.1.5, Theorem A]{Goresky1989}. 
They used stratified Morse theory to show that there is a finer gradation of the cohomolgoy groups of the complement indexed by intersections. 
In particular, they gave a closed form formula for these groups in terms of the homology groups of lower intervals in $L_{\A}$:
\begin{align}
\widetilde{H}^i(M_{\A}, \Z) \cong \bigoplus_X \widetilde{H}_{\mathrm{codim}(X)-2-i}(\Delta(\hat{0}, X), \Z).\label{gmf}
\end{align}
Here $X$ varies on all positive codimension intersections and on the right hand side are reduced simplicial homology groups of order complexes. 
We refer to the above expression as the Goresky-MacPherson formula for subspace arrangements. 
This result refines a theorem of Brieskorn and Orlik-Solomon concerning cohomology of hyperplane complements \cite[Lemma 5.91]{orlik92}.

The work of Goresky and MacPherson created a lot of activity; along with applications several authors gave alternate proofs and strengthening.
For an extensive survey of this subject (mainly for the work done till the 90's) we refer the reader to the papers of Bjorner \cite{Bjorner1994} and Vassiliev \cite{vas01}.
In the last decade the main focus has been the multiplicative structure of $H^{\ast}(M_{\A}, \Z)$. 
In particular, the arrangement of coordinate subspaces has applications in toric topology. 
See for example \cite{delong01b} and \cite[Chapter 8]{buch_panov02}.

We should mention the work of Vassiliev \cite{Vassiliev1992} and Ziegler-\u{Z}ivaljevi\'c \cite{zz93} where they give an alternate proof of the formula. 
Observe that the one-point compactification $\widehat{V_{\A}}$ of the link and $M_{\A}$ are complementary subspaces of $\widehat{\R^n} = S^n$.
In both these papers the homotopy type of $\widehat{V_{\A}}$ is expressed in terms of the combinatorial information. 
The Goresky-MacPherson formula then follows from Alexander duality.
The homotopy type of $\widehat{V_{\A}}$ was determined by Vassiliev using simplicial resolutions and by Ziegler-\u{Z}ivaljevi\'c using homotopy colimits.

The study of toric arrangements received prominence because of the work of De Concini and Procesi, \cite{copr05, cpbook09}. 
A toric arrangement is a finite collection of codimension-$1$ subtori in $(\C^*)^l$.
Besides being a natural generalization of hyperplane arrangements they have applications in diverse fields like algebraic geometry, combinatorics, partition functions and box splines. 
The last decade has seen a lot of activity in this field. 
We refer the reader to \cite[Section 1]{callegaro2017} for an excellent account of the state of the art in this area.
However, we do mention some of the relevant results concerning cohomology of the complement of a toric arrangement. 
The additive structure of cohomology was determined by Looijenga \cite{looi93} and also by De Concini and Procesi \cite{copr05}. 
It was later shown by D'Antonio and Delucchi in \cite{DAntonio2015}, using discrete Morse theory, that for the class of complexified toric arrangements the cohomology is torsion free.
Later the result was extended by Callegaro and Delucchi in \cite{callegaro2017} and now it is known that in all cases the cohomology groups are torsion free.
Recently, a complete multiplicative structure of the cohomology ring was obtained in \cite{cddmp18}.
See \cite{pagaria18} for examples of toric arrangements with isomorphic intersection data but non-isomorphic cohomlogy rings.
The rational formality (in the sense of Sullivan) of the complement was established by Dupont in \cite{Dupont2016}. 

The main aim of this paper is to establish a Goresky-MacPherson type formula in the context of toric arrangements. 
We consider finite collections of toric subspaces (not necessarily of codimension $1$). 
We show that there is a finer gradation of the cohomology groups of the complement. 
Moreover, this gradation is indexed by the underlying combinatorial data (Theorem \ref{mainthm}).
In order to prove this we find a closed form formula for the homotopy type of the one-point compactification of the link (Theorem \ref{thm03}). 
However, in the toric context we consider the `modified' link; the union of all the toric subspaces and the coordinate axes. 
The complement $M_{\A}$ and the one-point compactification of the modified link are complementary subspaces in $S^{2l}$.
So the formula is now a direct consequence of Alexander duality.
Our proof is based on the Ziegler-\v{Z}ivaljevi\'c approach for subspace arrangements.  
Towards the end we also make use of the spectral sequence developed by Petersen in \cite{petersen17} to get decomposition of the compactly supported cohomology (Theorem \ref{thmcpt}).
We pay special attention to arrangements of toric hyperplanes (Theorem \ref{toricgm}).
There are two immediate consequences of our result; first there is no torision in the cohomology (Corollary \ref{cor1}) and second is a theorem due to De Concini and Procesi (Corollary \ref{cor0}).
We end the paper with some discussion regarding future research directions. 

\section{Preliminaries}

In this section we sketch the proof, due to Ziegler-\v{Z}ivaljevi\'c \cite{zz93}, of the Goresky-MacPherson formula.
We begin by a brief review of topological aspects of posets; the reader's familiarity with posets and lattices is assumed.
In particular, we will make frequent use of the order complex of a poset so we state relevant definitions and facts here, one can find more details in \cite{wachs07}.

\begin{defn}\label{def5}
The M\"obius function of a poset $P$ is defined recursively on closed intervals of $P$ as follows:
\[\ds\mu(\alpha, \beta) = \begin{cases} 1, &\hbox{~if~} \alpha = \beta, \\
-\sum_{\alpha\leq \gamma< \beta}\mu(\alpha, \gamma), &\hbox{~if~} \alpha < \beta,\\
0 &\hbox{~if~}\alpha >\beta.\end{cases} \]
\end{defn}

The \emph{order complex} $\Delta(P)$ is an (abstract) simplicial complex whose $k$-simplices correspond to chain of length $k$. 
In this paper, we do not distinguish between abstract simplicial complexes and their geometric realizations.
The reduced homology of $P$, $\widetilde{H}_i(P)$, is defined as the reduced simplicial homology of $\Delta(P)$ with integer coefficients.
The cohomology of $P$ is defined analogously. 
We will often consider the order complex of an open interval $(\alpha, \beta) := \{\gamma\in P\mid \alpha< \gamma < \beta \}$. 
We denote its homology by $\redh_i(\alpha, \beta)$.

It is a basic fact, due to P. Hall, that the M\"obius function $\mu(\alpha, \beta)$ is the reduced Euler characteristic of $\Delta(\alpha, \beta)$:
\[\mu(\alpha, \beta) = \sum_i (-1)^i\; \mathrm{rank} \redh_i(\alpha, \beta). \]

A class of lattices called \emph{geometric lattices} plays an important role in this article. 
The intersection lattice of a central hyperplane arrangement (more generally, any proper interval in the intersection poset of a non-central arrangement) is a geometric lattice. 
A fundamental result due to Folkman states that for a geometric lattice $L$, the reduced homology groups $\redh_i(\hat{0}, \hat{1})$ vanish in all but the top dimension (i.e., dimension equal to $\mathrm{rank}(L) - 2$). 
In that dimension the homology is free of rank $|\mu(\hat{0}, \hat{1})|$. 
In fact, the homotopy type is that of a wedge of spheres of dimension $\mathrm{rank}(L) - 2$ (see \cite[Section 4.5]{orlik92} and \cite[Section 3.2.3]{wachs07}).

As introduced earlier we denote by 
$$\A := \{A_1,\dots, A_k\}$$ 
an arrangement of subspaces in $\R^n$. 
The \emph{intersection semilattice} $L_{\A}$ is the collection of all non-empty intersections of members of $\A$. 
The elements of $L_{\A}$ will be denoted by small Greek letters $\alpha, \beta$ and the corresponding intersections by $A_{\alpha}, A_{\beta}$ respectively.
The intersection semilattice is ordered by reverse inclusion: $\alpha\leq \beta \iff A_{\alpha}\supset A_{\beta}$.
This is a meet-semilattice with the least element $\hat{0} = \R^n$. 
In general $L_{\A}$ is not even graded.
The dimension function $d:L_{\A}\to \Z_{\geq 0}$ assigns to each $\beta$ the dimension of $A_{\beta}$.

The top element $\hat{1}$ exists if and only if the arrangement is central, i.e, the intersection of all the subspaces in $\A$ is non-empty. 
In such a situation the joins are also defined and $L_{\A}$ then becomes a lattice. 
If the context is clear we write $L$ instead of $L_{\A}$. 
Finally, use the notation $\lbu$ for $L_{\A}\setminus\{\hat{0}\}$.
We now (re)state the Goresky-MacPherson formula \cite[III.1.5]{Goresky1989}. 

\bt{gmformula01}
Let $\A$ be a subspace arrangement in $\R^n$, $M_{\A}$ be its complement. 
Then for every $i\geq 0$:
\[\widetilde{H}^i(M_{\A}) = \bigoplus_{\beta\in \lbu}\widetilde{H}_{n - d(\beta)-2-i}(\hat{0}, \beta). \]
\et

\be{gmex01}
Suppose that $\A$ is an arrangement of hyperplanes in $\R^n$.
In this case each closed interval $[\hat{0}, \beta]$ is a geometric lattice of rank $n-d(\beta)$.
The homology groups on the right are nonzero in the dimension $n-d(\beta)-2$, which implies $i=0$.
Consequently, one recovers the celebrated theorem of Zaslavsky \cite[Chapter 2]{orlik92}.
\ee

\be{gmex02}
Suppose $\A$ is a complex hyperplane arrangement in $\C^n\cong \R^{2n}$.
Again each closed interval is a geometric lattice of rank $n-d(\beta)$. 
And we get the following reformulation of Brieskorn's Lemma \cite[Lemma 5.91]{orlik92}. 
For $0\leq i\leq n$
\[\redh^{i}(\ma) \cong \bigoplus_{d(\beta) = n-i}\redh_{i-2}(\hat{0},\beta). \]
\ee 
The first proof of the Goresky-MacPherson formula was given as an application of stratified Morse theory. 
Some later reproofs involved analyzing, first, the topology of the link, which has more combinatorial structure, and then passing to the cohomology of the complement via Alexander duality.

We present a sketch of these ideas since the proof of our main result is based them. 
We begin by recalling some relevant notions from homotopy theory.
The interested reader can refer to \cite[Chapter 15]{kozlov_book}, \cite[Apendix]{zz93} and \cite{wzz99}.

We view a poset $P$ as a category in the standard way; elements of $P$ being objects and non-identity morphisms pointing down. 
Given a poset $P$ a $P$-diagram $\mathcal{D}$ of topological spaces is a covariant functor from $P\to\mathbf{Top}$, i.e., an assignment of spaces $\mathcal{D}(p)$ for every element $p\in P$ and continuous maps $\mathcal{D}(p>p'):\mathcal{D}(p)\to \mathcal{D}(p')$ for every order relation $p>p'$.
The colimit of a diagram $\mathcal{D}:P\to\mathbf{Top}$ is the topological space $\mathrm{colim}\mathcal{D}$ constructed from the disjoint union of all of $\mathcal{D}(p)$'s by identification of $x\in\mathcal{D}(p)$ with $\mathcal{D}(p>p')(x)\in \mathcal{D}(p')$, for all $p>p'$ and all $x\in\mathcal{D}(p)$.

Let $\mathcal{E}:P\to\mathbf{Top}$ be another $P$-diagram of topological spaces. 
A \emph{morphism of $P$-diagrams} $F:\mathcal{D}\to\mathcal{E}$ is a collection of continuous maps $F_p:\mathcal{D}(p)\to\mathcal{E}(p)$ indexed by $p\in P$ such that $\mathcal{E}(p>q)\circ F_p = F_q\circ \mathcal{D}(p>q)$.

Denote by $\mathrm{sd}(P)$ the poset of chains in $P$ and by $P_{\leq x}$ the subposet of elements below $x$.

\begin{defn}\label{hocolim}
The \emph{homotopy colimit} of a diagram $\mathcal{D}: P\to \mathbf{Top}$, denoted by $\mathrm{hocolim}(\mathcal{D})$, is the colimit of the functor $\Delta(\mathcal{D}):\mathrm{sd}(P)\to \mathbf{Top}$ defined by:
\begin{itemize}
\item on the elements: $\Delta(\mathcal{D})(x_1>\cdots>x_i) = \Delta(P_{\leq x_i})\times \mathcal{D}(x_1)$;
\item on the morphisms:\[\Delta(\mathcal{D})((x_1>\cdots>x_i)\to(x_{j_1}>\cdots> x_{j_k})) = \mathit{I}(x_i\hookrightarrow x_{j_k})\times \mathcal{D}(x_1\to x_{j_1}). \]
\end{itemize}
\end{defn}

Note that a morphism of $P$-diagrams induces a continuous map between the corresponding homotopy colimits.
We are mainly interested in the following three properties of homotopy colimits. 
\begin{enumerate}
\item The \emph{homotopy lemma} \cite[Proposition 3.7]{wzz99}; it says that the homotopy type of the homotopy colimit doesn't change if spaces in a diagram are changed up to homotopy.
\begin{lemma}\label{htpy}
Let $F:\mathcal{D}\to\mathcal{E}$ be a morphism of two $P$-diagrams such that each $F_p$ is a homotopy equivalence.
Then the map induced by $F$ on the homotopy colimits is also a homotopy equivalence.
\end{lemma}
\item The \emph{projection lemma} \cite[Lemma 4.5]{wzz99}; identifies a situation in which the natural map from $\mathrm{hocolim}(\mathcal{D})\to \mathrm{colim}(\D)$ is a homotopy equivalence.
\begin{lemma}\label{proj}
Let $\mathcal{D}$ be a diagram indexed by a poset such that all the induced maps are inclusions that are also closed cofibrations. Then the collapsing map from $\mathrm{hocolim}(\mathcal{D})\to \mathrm{colim}(\D)$ is a homotopy equivalence.
\end{lemma}
\item The \emph{wedge lemma} \cite[Lemma 4.9]{wzz99}, we include the technical statement instead of its intuitive description.
\begin{lemma}\label{wedgel}
Suppose $P$ is a poset with maximal element $\hat{1}$. Let $\D$ be a $P$-diagram such that all the maps between spaces are constant. Then 
\[\mathrm{hocolim}(\D) \simeq \bigvee_{\alpha\in P} \left(\Delta(P_{<\alpha})\ast \D(\alpha) \right), \]
where ``$\ast$'' denotes the topological join of two spaces.
\end{lemma}
\end{enumerate}

\begin{defn}
Given a subspace arrangement $\A$ the \emph{subspace diagram} $\mathcal{D}_{\A}: \lbu\to \mathbf{Top}$ is defined as follows:
\begin{itemize}
\item an object $\beta$ is sent to the corresponding intersection $A_{\beta}$;
\item a morphism $\beta > \alpha$ is sent to the corresponding inclusion $A_{\beta} \hookrightarrow A_{\alpha}$.
\end{itemize}
\end{defn}

Now we can proceed with the proof of the Goresky-MacPherson formula. 
The first step is to determine the homotopy type of the link. 

\begin{thm}[{\cite[Theorem 2.1]{zz93}}]
For a subspace arrangement $\A$ the homotopy type of the link is given by 
\[V_{\A}\simeq \Delta(\lbu).\]
\end{thm}

In order to arrive at the above result, first note that because of the projection lemma $\mathrm{hocolim}(\D_{\A})$ is homotopy equivalent to ${V}_{\A}$. 
On the other hand the order complex of $\lbu$ is homotopy equivalent to $\mathrm{hocolim}(\mathcal{D}_{\A})$ because of the homotopy lemma.

The homotopy type of the compactification $\widehat{V}_{\A}$ is also combinatorially determined. 
Here $S^i$ denotes the $i$-sphere. 
Also, it is to be understood that $\emptyset\ast A = A$ for all spaces $A$.
\begin{thm}[{\cite[Theorem 2.2]{zz93}}]\label{linkthm}
For every subspace arrangement:
\[ \widehat{V}_{\A} \simeq \bigvee_{\beta\in \lbu} \left( \Delta(\hat{0}, \beta)\ast S^{d(\beta)}\right).\]
Consequently the homology is given by 
\[\widetilde{H}_i(\widehat{V}_{\A}) \cong \bigoplus_{\beta\in \lbu} \tilde{H}_{i-d(\beta)-1}(\hat{0}, \beta). \]
\end{thm}

The proof of the above theorem involves modifying the diagram $\mathcal{D}_{\A}$. 
Add an extra element $\hat{1}$ to $\lbu$ which will be the index of the point at infinity. 
Now to each element $\beta\in \lbu$ assign the sphere $S^{d(\beta)}$ which is homeomorphic to the one-point compactification $\widehat{A}_{\beta}$ of $A_{\beta}$. 
To $\hat{1}$ assign the point at infinity; it gets included in all the spheres. 
The maps in this diagram are null homotopic; in fact, the union $\cup_{\alpha>\beta}\widehat{A}_{\alpha}$ is contained in a closed disc in $\widehat{A}_{\beta}$.
Consequently one may replace this diagram with a new diagram consisting of spheres of appropriate dimensions and base point preserving constant maps. 
The result now follows from the application of the wedge lemma (Lemma \ref{wedgel} above). 

Note that the compactified link $\widehat{V}_{\A}$ and the complement $M_{\A}$ are complementary subspaces of $\widehat{\R^n}$. 
Hence the Goresky-MacPherson formula is a consequence of Alexander duality:
\[\widetilde{H}^i(M_{\A}) \cong \widetilde{H}_{n-i-1}(\widehat{V}_{\A}). \]

The approach taken by Vassiliev in \cite{Vassiliev1992} to determine the homotopy type of $\widehat{V}_{\A}$ is slightly different. 
Vassiliev constructed a `simplicial blow up' of the intersections of subspaces:
Take $N$ sufficiently large number such that each subspace $A_{\beta}$ can be embedded in a generic position. 
Now for every $x\in V_{\A}$ let $S(x)$ denote the convex hull of the images of $x$ in 
$\R^N$. 
The simplicial resolution $S(\A)$ of $V_{\A}$ is the union of all $S(x)$'s. 
He goes on to prove that $\widehat{S(\A)} \simeq \widehat{V}_{\A}$. 
Then he uses the stratified Morse theory to arrive at the closed form formula for $\widehat{S(\A)}$ (see also \cite[Section 8]{vas01}). 
It was later proved by Kozlov in \cite{Kozlov2002} that the simplicial resolution itself is a homotopy colimit. 
Moreover he constructed an explicit deformation retraction from $\widehat{S(\A)}$ to the homotopy colimit constructed by Ziegler-\v{Z}ivaljevi\'c.

\section{Toric arrangements}\label{sec2}

We now turn to the main object of this paper; toric arrangements. 
The standard $l$-dimensional complex torus is the space $(\C^*)^l$ of $l$-tuples of nonzero complex numbers. 
The torus is a group under multiplication of coordinates. 
In fact, it is an affine algebraic variety with the ring of \emph{Laurent polynomials} $\C[z_1^{\pm 1},\dots, z_l^{\pm 1}]$ as its coordinate ring. 

\bd{def1}
The \textit{character} of a torus is a multiplicative homomorphism $\chi\colon (\C^*)^l\to \C^*$ given by the evaluation of Laurent monomials  
\[\chi(z_1, \dots, z_l) = z_1^{n_1}\cdots z_l^{n_l}, \quad n_i\in \Z, \forall i. \] \ed

The following are some well-known facts regarding tori (see \cite[Section 5.2]{cpbook09}). 
The set $\Lambda$ of characters of $(\C^*)^l$ is a free abelian group isomorphic to $\Z^l$. 
Conversely, for a finitely generated abelian group $\Lambda$ of rank $l$ the variety $T^l_{\C} := \mathrm{Hom}(\Lambda, \C^*)$ is isomorphic to the product of an $l$-torus and a finite abelian group isomorphic to the torsion subgroup of $\Lambda$.

The connected, topologically closed subgroups (called the toric subgroups) of $\ctor{l}$ are isomorphic to $k$-tori for some $k\leq l$. 
In general a closed subgroup  of $\ctor{l}$ is isomorphic to $\ctor{k}\times A$ where $A$ is a finite abelian group. 
Any coset of a toric subgroup is homeomorphic to the group; thus topologically it is a torus. 
Let $W$ be a closed subgroup of $\ctor{l}$ and $W_0$ be the connected component containing identity. 
Then $W_0$ is a toric subgroup, the quotient $W/W_0$ is a finite abelian subgroup and $W$ is the union of the $|W/W_0|$ distinct cosets. \par 

There is a one-to-one correspondence between closed subgroups of $\ctor{l}$ and subgroups of $\Lambda$. 
Such subgroups are determined by integer matrices. An integer matrix $A$ of size $m\times l$ determines a mapping from $\ctor{l}$ to $\ctor{m}$. 
The kernel $W$ of this mapping is a closed subgroup of $\ctor{l}$ and every subgroup arises in this manner. 
Consequently $W$ depends only on the subgroup (i.e., the sub-lattice) of $\Lambda$ generated by the rows of $A$. 
Without loss of generality one can assume that rows of $A$ furnish a basis for this sub-lattice so that $m\leq l$.
By a \emph{toric subspace} we mean a translate of a closed subgroup.

Because of their importance, we first focus on collections consisting of codimension-$1$ subtori.
Given a character $\chi\in\Lambda$ and a non-zero complex number $c$ a \textit{toric hypersurface} $K_{\chi, c}$ is defined as the level set of $\chi$, i.e., $K_{\chi, c} := \{z\in (\C^*)^l\mid \chi(z) = c \}$. 
A toric hypersurface is a translate of a toric subgroup of codimension-$1$.

\bd{def2}
A \emph{toric arrangement} is a finite collection 
\[\A := \{K_{\chi_1, c_1}, \dots,  K_{\chi_n, c_n}\mid \chi_i \in \Lambda, c_i\in \C^* ~\hbox{for}~ 1\leq i\leq n\}\]
of toric hypersurfaces in $(\C^*)^l$.
\ed

For notational simplicity we write $K_i$ instead of $K_{\chi_i, c_i}$. 
Without loss of generality we assume that each toric hypersurface in $\A$ is connected, i.e., each character is primitive (recall that a character is primitive if all the exponents are relatively prime). 


\begin{defn}\label{intersection poset}
The \textit{intersection poset} $L_{\A}$ (or the \emph{poset of layers}) is the set of all connected components of all nonempty intersections of the toric hypersurfaces in $\A$ ordered by reverse inclusion.
The elements of $L_{\A}$ are called \textit{components of the arrangement}.
\end{defn}

\begin{figure}[ht!]
\begin{center}
\begin{tikzpicture}[scale=3.5]
\draw[-] (0,0) rectangle ++(1,1);
\draw[-] (0,1) -- (0.5,0);
\draw[-] (0,1) -- (1,0.5);
\draw[-] (1,0) -- (0.5,1);
\draw[-] (1,0) -- (0,0.5);
\filldraw[-] (1/3,1/3) circle (0.5pt);
\filldraw[-] (2/3,2/3) circle (0.5pt);
\filldraw[-] (1,0) circle (0.5pt);
\node at (1.07,-0.07) {$p_1$};
\node at (0.4,0.4) {$p_2$};
\node at (0.73,0.73) {$p_3$};
\node at (0.5,-0.07) {$K_2$};
\node at (0.5,1.07) {$K_2$};
\node at (-0.07,0.5) {$K_1$};
\node at (1.12,0.5) {$K_1$};
\node at (0.5,-0.3) {$\A$};
\end{tikzpicture}
\quad\quad\quad
\begin{tikzpicture}[scale=3]
\filldraw[-] (1/2,0) circle (0.6pt);
\filldraw[-] (0,1/2) circle (0.6pt);
\filldraw[-] (0,1) circle (0.6pt);
\filldraw[-] (1/2,1) circle (0.6pt);
\filldraw[-] (1,1) circle (0.6pt);
\filldraw[-] (1,1/2) circle (0.6pt);
\draw[-] (1/2,0) -- (0,1/2);
\draw[-] (0,1/2) -- (0,1);
\draw[-] (0,1/2) -- (1/2,1);
\draw[-] (0,1/2) -- (1,1);
\draw[-] (1/2,0) -- (1,1/2);
\draw[-] (1,1/2) -- (1,1);
\draw[-] (1,1/2) -- (1/2,1);
\draw[-] (1,1/2) -- (0,1);
\node at (0.4,-0.1) {$(\C^*)^2$};
\node at (-0.12,1/2) {$K_1$};
\node at (1.12,1/2) {$K_2$};
\node at (0,1.07) {$p_1$};
\node at (1/2,1.07) {$p_2$};
\node at (1,1.07) {$p_3$};
\node at (0.5,-0.36) {$L_{\A}$};
\end{tikzpicture}
\end{center}
\caption{The real part of a toric arrangement in $(\C^*)^2$ and its intersection poset.}
\label{fig1}
\end{figure}
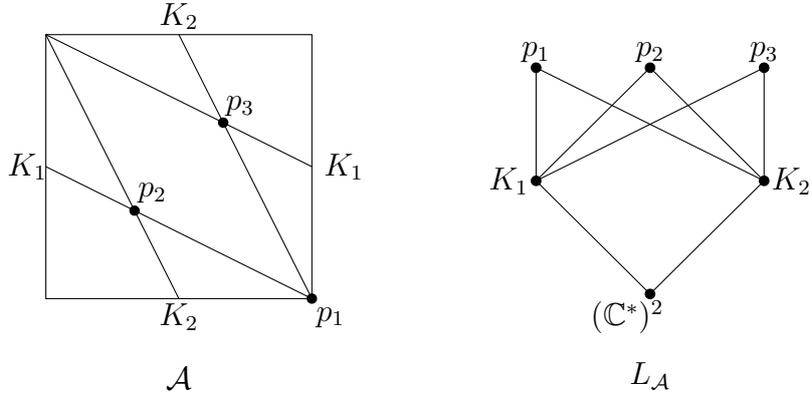

\begin{figure}[ht!]
\begin{center}
\begin{tikzpicture}[scale=3.5]
\draw[-] (0,0) rectangle ++(1,1);
\draw[-] (0,1) -- (0.5,0);
\draw[-] (0,1) -- (1,0.5);
\draw[-] (1,0) -- (0.5,1);
\draw[-] (1,0) -- (0,0.5);
\draw[-] (0,0) -- (1,1);
\filldraw[-] (1/3,1/3) circle (0.5pt);
\filldraw[-] (2/3,2/3) circle (0.5pt);
\filldraw[-] (1,0) circle (0.5pt);
\node at (1.07,-0.07) {$p_1$};
\node at (0.45,0.35) {$p_2$};
\node at (0.78,0.68) {$p_3$};
\node at (0.5,-0.07) {$K_2$};
\node at (0.5,1.07) {$K_2$};
\node at (-0.07,0.5) {$K_1$};
\node at (1.09,0.5) {$K_1$};
\node at (1.05,1.05) {$K_3$};
\node at (0.5,-0.3) {$\A$};
\end{tikzpicture}
\quad\quad\quad
\begin{tikzpicture}[scale=3]
\filldraw[-] (1/2,0) circle (0.6pt);
\filldraw[-] (0,1/2) circle (0.6pt);
\filldraw[-] (1/2,1/2) circle (0.6pt);
\filldraw[-] (0,1) circle (0.6pt);
\filldraw[-] (1/2,1) circle (0.6pt);
\filldraw[-] (1,1) circle (0.6pt);
\filldraw[-] (1,1/2) circle (0.6pt);
\draw[-] (1/2,0) -- (0,1/2);
\draw[-] (1/2,0) -- (1/2,1/2);
\draw[-] (1/2,0) -- (1,1/2);
\draw[-] (0,1/2) -- (0,1);
\draw[-] (0,1/2) -- (1/2,1);
\draw[-] (0,1/2) -- (1,1);
\draw[-] (1,1/2) -- (1,1);
\draw[-] (1,1/2) -- (1/2,1);
\draw[-] (1,1/2) -- (0,1);
\draw[-] (1/2,1/2) -- (1,1);
\draw[-] (1/2,1/2) -- (1/2,1);
\draw[-] (1/2,1/2) -- (0,1);
\node at (0.4,-0.1) {$(\C^*)^2$};
\node at (-0.12,1/2) {$K_1$};
\node at (0.38,1/2) {$K_2$};
\node at (1.12,1/2) {$K_3$};
\node at (0,1.07) {$p_1$};
\node at (1/2,1.07) {$p_2$};
\node at (1,1.07) {$p_3$};
\node at (0.5,-0.36) {$L_{\A}$};
\end{tikzpicture}
\end{center}
\caption{The real part of a toric arrangement and the intersection poset.}
\label{fig2}
\end{figure}
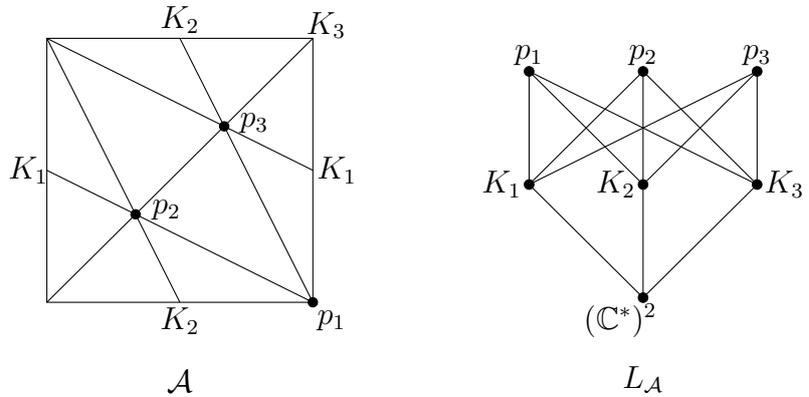

The intersection poset is in general not a semilattice.
However, it is certainly a graded poset; the rank of every element is the codimension of the corresponding intersection.
The ambient torus $(\C^*)^l$ is the least element $\hat{0}$. 
As before, the elements of $L_{\A}$ are denoted by $\alpha, \beta$ and the corresponding intersections by $W_{\alpha}, W_{\beta}$ respectively. 
We also have the (complex) dimension function $d: L_{\A}\to \Z_{\geq 0}$.
Note that for every $\beta\in L_{A}$ the intersection $W_{\beta}$ isomorphic to $(\C^*)^{d(\beta)}$.
Every proper interval in $L_{\A}$ is a geometric lattice, see \cite[Theorem 3.11]{ehr09}.
A few examples before we proceed. 

\be{ex1} An arrangement in $\C^*$ is just a collection of finitely many points. \ee

\be{ex2}
Let $\A = \{z_1z_2^2 =1, z_1^2z_2=1 \}$ be a toric arrangement in $(\C^*)^2$. 
The toric hypersurfaces $K_1$ and $K_2$ intersect in three points $p_1 = (1, 1), p_2 = (\omega, \omega)$ and $p_3 = (\omega^2, \omega^2)$.
Where $\omega$ is a cube root of unity.
Figure \ref{fig1} shows the ``real part'' of the arrangement, i.e., intersection with the compact torus $(S^1)^2\subset (\C^*)^2$. 
It is visualized as a quotient of the unit square.  
Along with it is the Hasse diagram of $L_{\A}$.
\ee

\be{ex3}
Now consider the arrangement formed by adding the toric hypersurface $K_3$ given by equation $z_1z_2^{-1} = 1$ to the previous arrangement. 
The hypersurfaces $K_1, K_2$ and $K_3$ intersect in the same three points $p_1, p_2, p_3$ as above.
Figure \ref{fig2} shows the ``real part'' of the arrangement and the associated intersection poset. 
\ee

\be{ex4} The braid arrangement in $(\C^*)^l$ is the collection of $\binom{l}{2}$ toric hypersurfaces given by the following equations 
\[\A = \{z_i z_j ^{-1} = 1 \mid 1\leq i < j \leq l\}. \]
The intersection poset in this case is the partition lattice $\Pi_n$. 
\ee

\bd{def3}
The \textit{complement} of a toric arrangement $\A$ in $(\C^*)^l$ is
\[\ma := (\C^*)^l \setminus \bigcup_{i=1}^{n} K_i.\]
\ed
For example, the complement of a toric arrangement in $\C^*$ has the homotopy type of wedge of circles.
The complement of the braid arrangement is the configuration space of $l$ ordered points in $\C^*$.

Let $W_{\beta}$ be a component of a toric arrangement $\A$. 
For any $z\in W_{\beta}$ let $\A^z[\beta]$ denote the hyperplane arrangement in a coordinate neighborhood of $z$. 
More intrinsically, $\A^z[\beta] = \{T_z(K)\mid W_{\beta}\subseteq K \}$ in the tangent space $T_z(\C^*)^l \cong \C^l$. 
Note that if $z$ is a point such that it does not belong to any other component contained in $W$ then the intersection data of $\A^z[\beta]$ is independent of the choice of $z$. 
One can safely disregard the reference to $z$ as we are mainly interested in the combinatorics. 
We refer to $\A[\beta]$ as the \textit{local arrangement} at $W_{\beta}$. 
The complement $M(\A_{\beta})$ is the \textit{local complement} at $W_{\beta}$.
Finally, for $0\leq j\leq l$ define numbers $N_j(\A)$ as follows:
\[ N_j(\A) := \sum_{\substack{\beta\in L_{\A}\\ d(\beta) = l-j}} \rank H^j(M(\A[\beta])),\]
i.e., the sum of $j$-th Betti numbers of all local complements at codimension $j$ components of $\A$. 
The reader familiar with hyperplane arrangements will realize that the summands are precisely cardinality $j$ \emph{no-broken-circuits} of local arrangement $\A[\beta]$. 
Further details can be found in \cite[Section 3.1]{orlik92} and in \cite[Section 3.4]{copr05} (in the context of toric arrangements).

The cohomology groups of $M_{\A}$ were computed and by De Concini and Procesi in \cite[Theorem 4.2]{copr05} using an algebraic De Rham complex. 
We state an equivalent version of their theorem below.
\bt{thm0} Given a toric arrangement $\A$ the Poincar\'e polynomial of its complement is given by
\[\ds \mathrm{Poin}(\ma, t) = \sum_{j=0}^l N_j (t+1)^{l-j} t^j.\]
In particular, the $j$-th Betti number is $\sum_{i=0}^j N_i \binom{l-i}{j-i}$.
\et 

\be{dpex}
Let us apply the above theorem to the arrangement in Example \ref{ex4}. 
Observe that the local arrangement at every point is the arrangement of $3$ concurrent lines.
Whereas the local arrangement at every $K_i$ is isomorphic to the arrangement of a point in $\C$. 
Hence, $N_0 = 1, N_1 = 3$ and $N_2 = 6$. 
Plugging these numbers in the above formula we get 
$\rank H^0(\ma) = 1, \rank H^1(\ma) = 5, \rank H^2(\ma) = {10} $.
\ee 

We now proceed to prove the analogue of the Goresky-MacPherson formula in this context. 
First, a simple observation. 
Given a toric arrangement $\A$ in $(\C^*)^l$ we have
\[ \ma \cong \C^l \setminus \left(\{z_i = 0\mid 1\leq i\leq l\} \cup \bigcup_{i=1}^n K_i\right).\]
The observation motivates the following definition.

\bd{tlink}
Given a toric arrangement $\A$ in $(\C^*)^l$ the \emph{(modified) link} $\ua$ is defined as 
\[ \ua := \{z_i = 0\mid 1\leq i\leq l\} \cup \bigcup_{i=1}^n K_i.\]
\ed
Forgetting the complex structure on $\C^l$ and passing on to the compactification we get 
\begin{align}
\ma \cong \widehat{\R^{2l}}\setminus \widehat{\ua}.
\end{align}
Hence our strategy is to express the homotopy type of $\widehat{\ua}$ in terms of the underlying combinatorics and then use the Alexander duality to arrive at the Goresky-MacPherson formula. 
We begin by proving a couple of lemmas. 

\bl{lem01}
Let $B_l$ denote the union of all the coordinate axes in $\C^l$. 
Then 
\[\widehat{B_l}\simeq \bigvee_{i=1}^l\left( \bigvee_{\binom{l}{i}} S^{2l-i-1} \right) \]
\el
\bpr 
Notice that $B_l$ is the link of the arrangement of coordinate hyperplanes in $\C^l$. 
The associated intersection poset is the Boolean lattice, which we denote by $2^{[l]}$.
Each lower interval in $2^{[l]}$ is homotopic to the sphere of appropriate dimension.
The result follows from Theorem \ref{linkthm} as join of two spheres is again a sphere.
\epr

\bl{lem02}
The homotopy type of the one-point compactification $\widehat{(\C^*)^d}$ has the following closed form
\[ \widehat{(\C^*)^d} \simeq \bigvee_{i=0}^d \left(\bigvee_{\binom{d}{i}}S^{d+i} \right).\]
\el 
\bpr 
Recall that for two locally compact, Hausdorff spaces $A, B$ we have the following homeomorphism
\begin{align}
\widehat{A\times B}\cong \widehat{A}\wedge\widehat{B} \label{cpth}
\end{align}
that maps the compactification point to the smash point and is identity on $A\times B$.

Observe that $\widehat{\C^*}$ is homotopic to $S^1\vee S^2$. 
To see this, first consider $\C^*$ as a subset of the Riemann sphere; i.e., the sphere with both poles missing.
Now add the missing two points and and then identify them. 
As the smash product distributes over the wedge, the following homotopy equivalence establishes the result
\begin{align*}
\widehat{(\C^*)^d} &\simeq \left( S^1\vee S^2\right)^{\wedge d}\\
				   &= \bigvee_{i_k\in\{1, 2\}} S^{i_1}\wedge\cdots\wedge S^{i_d}\\
                   &= \bigvee_{i=0}^d\bigvee_{\binom{d}{i}} S^{d+i}.
\end{align*}
\epr 

\br{rem01}
For notational simplicity we denote by $\tdhat{d}$ the wedge of spheres on the right hand side of the above formula. 
Also, note that $\tdhat{0} = S^0$.
\er
Given a toric arrangement $\A$ in $(\C^*)^l$ the \emph{modified intersection poset} is defined as 
\[ P_{\A} := \lbu \amalg 2^{[l]}\setminus\{\emptyset \},\]
here $2^{[l]}$ denotes the Boolean lattice. 
\bd{toricd}
The \emph{toric subspace diagram} $\D_{\A}: P_{\A}\to\mathbf{Top}$ is defined as follows:
\begin{enumerate}
\item If $\beta\in\lbu$ then $\D_{\A}(\beta)$ is $W_{\beta}$, the corresponding intersection. 
If $\beta\in 2^{[l]}\setminus\{\emptyset \}$ then  $\D_{\A}(\beta)$ is $\{z_i=0\mid i\in \beta \}$. 
\item The maps are the inclusions of subspaces.
\end{enumerate}
\ed

\bl{proj02}
The homotopy colimit of $\D_{\A}$ is homotopy equivalent to the modified link $\ua$.
\el

\bpr 
Follows from the projection lemma.
\epr 



\begin{proposition}\label{cnbd}
The point $\{\infty\}$ has a compact cone neighborhood in $\widehat{(\C^{\ast})^d}$.
\end{proposition}
\bpr 
For $n\geq 1$ define
\[U_n := \{(z_1,\dots,z_d)\in (\C^{\ast})^d \mid |z_i|> n,\, \forall i \}. \]
Each $U_n$ is a neighborhood of $\{\infty\}$ in $\widehat{(\C^{\ast})^d}$ and there are inclusions: $U_1\supset U_2\supset \cdots$.
For every $n$ there is a deformation retraction of $\overline{U_n}$ on to $\overline{U_{n+1}}$ that fixes $\overline{U_{n+2}}$. 
Paste these retractions to contract $\overline{U_1}$ on to $\{\infty\}$ in $\widehat{(\C^{\ast})^d}$.
\epr 
\begin{remark}
In fact, the point $\{\infty\}$ has arbitrarily small contractible neighborhoods which are wedge of $2d$ cones with $\{\infty\}$ as the common cone point.
For example, let 
$$C_n := \{z\in (\C^{\ast})^d \mid 1/n < |z_i| < n+1,\, \forall i\}.$$ 
The closures $\overline{C_n}$'s form an exhaustion by compact sets and the complement of each them has $2d$ connected components (i.e., the $d$-torus has $2d$ ends). 
By similar arguments as above one can construct a cone over $2d$ components of $\partial \overline{C_n}$. 
\end{remark}

\begin{proposition}\label{incone}
Suppose $W_{\beta}\subsetneq (\C^{\ast})^d$ is a toric subspace. 
Then its one point compactification $\widehat{W_{\beta}}$ is contained in a compact cone neighborhood of $\{\infty\}$ in $\widehat{(\C^{\ast})^d}$.
\end{proposition}

\bpr 
Recall that a toric subspace can be specified by a set of characters.
Assume that the given ${W_{\beta}}$ is the intersection of $m$ toric hyperplanes whose defining characters are $z_1^{n_{j1}}\cdots z_1^{n_{jd}} = c_j$ for $j = 1,\dots, m$.

Define  
\[\mathcal{R}(\beta):= \{x\in \R^d_{>0}\mid x_1^{n_{j1}}\cdots x_d^{n_{jd}} - |c_j| = 0,\; j=1,\dots, m\}. \]
Since $\mathcal{R}(\beta)$ is closed we can choose two $d$-tuples $(r_1,\dots, r_d)$ and $(R_1,\dots, R_d)$ such that they are in the same path component of $\R^d_{>0}\setminus \mathcal{R}(\beta)$, $r_i < R_i$ and both are rationals for every $i$.

Consider the following compact subset of $(\C^{\ast})^d$:
\[\mathcal{K}(\beta) := \{z\in (\C^{\ast})^d \mid r_i \leq |z_i| \leq R_i,\; i=1,\dots, d\} \]
(even though the definition depends on the choice of $r_i$'s and $R_i$'s we choose to ignore them in the notation).
The construction implies that $ \mathcal{K}(\beta)\cap W_{\beta} = \emptyset$.
Hence $W_{\beta}$ is contained in a neighborhood of an end of $(\C^{\ast})^d$. 
Equivalently, the compactification $\widehat{W_{\beta}}$ is contained in a compact cone on $\partial\mathcal{K}(\beta)$ with $\{\infty\}$ as its cone point. 
\epr 
We now determine a closed form formula for the homotopy type of the one point compactification of the modified link. 
\bt{thm03}
The homotopy type of the compactified link of $\A$ is given by
\begin{align}
\widehat{\ua} \simeq \widehat{B_l}\vee \bigvee_{\beta\in\lbu}\left[ \bigvee_{j=0}^{d(\beta)}\left(\bigvee_{\binom{d(\beta)}{j}}\Sigma^{d(\beta)+j+1} \Delta(\hat{0}, \beta)\right) \right]\label{clink}
\end{align}
where $\Sigma^i$ denotes the $i$-fold reduced suspension. 
\et
\bpr 
The proof is similar to the one for subspace arrangements (Theorem \ref{linkthm}) originally given by Ziegler-\v{Z}ivaljevi\'c in \cite{zz93}.

Construct a new diagram such that the spaces indexed by $P_{\A}$ are replaced by their one-point compactifications. 
In particular, we have $\D': P_{\A}\amalg\{\hat{1}\}\to \mathbf{Top}$ given by setting $\D'(\beta) = \widehat{\D_{\A}(\beta)}$ for $\beta\in P_{\A}$ and $\D'(\hat{1}) = \{\infty \}$.
The maps between the spaces are inclusions.

We claim that these inclusions are null homotopic. 
This is clear for spaces indexed by $2^{[l]}\setminus\emptyset$.
For spaces indexed by $\lbu$ it follows from Proposition \ref{incone}.
By projection lemma $\mathrm{hocolim}(\D')\simeq\widehat{\ua}$.

Now construct one more diagram, say $\mathcal{E}$, indexed by $P_{\A}\amalg\{\hat{1}\}$: if $\beta\in \lbu$ then $\mathcal{E}(\beta)=\tdhat{d(\beta)}$ otherwise ${\beta}\in 2^{[l]}\setminus\{\emptyset\}$ and $\mathcal{E}(\beta)$ is the wedge of spheres homotopic to the one-point compactification of  $\{z_i = 0\mid i\in\beta\}$ as obtained in Lemma \ref{lem01}. 
Finally, send $\hat{1}$ to $\{\ast\}$. 
The maps in this new diagram are the constant maps to the base point.

We now construct a morphism $\mathbf{h}:\mathcal{D}'\to\mathcal{E}$ as follows. First, a claim: for a given $\alpha\in\lbu$ the union $\cup_{\beta>\alpha}{\D'_{\A}(\beta)}$ is contained in a compact cone neighborhood of $\{\infty\}$ in $\D'_{\alpha}$.
In order to prove this claim we employ the same strategy as above. 
Construct a compact set $\mathcal{K}(\alpha) \subset \D_{\A}(\alpha)$ by choosing two points in $\R^d_{>0}\setminus \cup \mathcal{R}(\beta)$ (here the union is over all $\beta > \alpha$).
The claim now follows since the union of compactifications is contained in the compact cone on $\partial\mathcal{K}(\alpha)$. 
Choose a homotopy equivalence $h_{\alpha}:\D'(\alpha)\to\mathcal{E}(\alpha)$ that contracts this compact cone and sends it to the base point of $\mathcal{E}(\alpha)$.

For $\alpha\in 2^{[l]}\setminus\{\emptyset\}$ choose a homotopy equivalence that contracts the closed disk containing the union of smaller-dimensional spheres to the base point. 
The homotopy lemma now implies that $\mathrm{hocolim}(\mathcal{E})$ is homotopy equivalent to $\widehat{\ua}$. 

The application of the wedge lemma to $\mathrm{hocolim}(\mathcal{E})$ gives us the following decomposition:
\begin{align*}
\widehat{\ua} \simeq \widehat{B_l}\vee \bigvee_{\beta\in\lbu}\left(\Delta(\hat{0}, \beta)\ast \tdhat{d(\beta)} \right).
\end{align*}
Recall that for well-pointed spaces the join distributes over the wedge. Hence, 
\begin{align*}
\Delta(\hat{0}, \beta)\ast \tdhat{d(\beta)} &\simeq 
\Delta(\hat{0}, \beta)\ast \bigvee_{j=0}^{d(\beta)}\bigvee_{\binom{d(\beta)}{j}} S^{d(\beta)+j} \\
				 &\simeq \bigvee_{j=0}^{d(\beta)}\bigvee_{\binom{d(\beta)}{j}} \Delta(\hat{0}, \beta)\ast S^{d(\beta)+j}
\end{align*}
The result then follows from the fact that join with an $i$-sphere is homotopic to the $(i+1)$-fold reduced suspension.
\epr 


\bc{cor1}
Let $\A$ be a toric arrangement and $\ma$ be the associated complement. 
Then the cohomology groups $H^i(\ma)$ are torsion free.
\ec

\bpr
Note that since each interval is a geometric lattice the corresponding order complex $\Delta(\hat{0}, \beta)$ is homotopic to a wedge of spheres \cite[Theorem 4.109]{orlik92}. 
Their dimension being $ l-d(\beta) -2$ and their number being $|\mu(\hat{0}, \beta)|$. 
Consequently we have
\[
\widehat{\ua} \simeq \widehat{B_l}\vee \bigvee_{\beta\in\lbu}\left[ \bigvee_{j=0}^{d(\beta)}\left(\bigvee_{\binom{d(\beta)}{j}+|\mu(\hat{0},\beta)|}\; S^{l+j-1}\right) \right].
\]
It means that the homotopy type of $\widehat{\ua}$ is again a wedge of spheres.
\epr

Finally, the toric analog of the Goresky-MacPherson formula.
\bt{toricgm}Let $\A$ be a toric arrangement in $(\C^*)^l$. 
Then we have the following decomposition of the cohomology groups of the complement $\ma$.
\begin{align}
\tld{H}^i(\ma) \cong 
\bigoplus_{\beta\in L_{\A}}\tld{H}_{l-d(\beta)-2}(\hat{0}, \beta)\otimes \tld{H}^{l-i}(W_{\beta})
\label{heq01}
\end{align}
\et
\bpr 
By Alexander duality we have $\tld{H}^i(\ma)\cong \tld{H}_{2l-1-i}(\widehat{\ua})$. 
The homology groups of $\widehat{\ua}$ are direct sum of homology of factors in the wedge decomposition in Equation \ref{clink}. 
In particular we have, 
\begin{align}
\tld{H}_j(\widehat{\ua}) \cong \redh_j(\widehat{B_l})\oplus \left(\bigoplus_{\beta\in\lbu} \redh_j(\Delta(\hat{0}, \beta)\ast \tdhat{d(\beta)})\right).\label{heq02}
\end{align}
Hence we determine each summand separately and plug it back in Equation \ref{heq02}.
First, observe that $\widehat{B_l}$ and $(\C^*)^l$ are complementary subspaces, hence 
\[\redh_{2l-1-i}(\widehat{B_l}) \cong \redh^i((\C^*)^l)\cong \redh^{l-i}(W_{\hat{0}}). \]
We now turn to the join of spaces, for notational simplicity let $k=2l-i-1$. 
\begin{align*}
\redh_k(\Delta(\hat{0}, \beta)\ast \tdhat{d(\beta)}) &\cong \redh_k\left(\bigvee_{j=0}^{d(\beta)}\bigvee_{\binom{d(\beta)}{j}} \Sigma^{d(\beta)+j+1} \Delta(\hat{0}, \beta)\right) \\
                 &\cong \bigoplus_{j=0}^{d(\beta)} \,\bigoplus_{\binom{d(\beta)}{j}} \redh_{k-d(\beta)-1-j}(\hat{0}, \beta).
\end{align*}
Due to Folkman's theorem homology of $\Delta(\hat{0}, \beta)$ is concentrated in degree $l-d(\beta)-2$. 
Hence the summand corresponding to $j = k-l-1$ survives and all others vanish. 
Hence,
\begin{align*}
\redh_{2l-i-1}(\Delta(\hat{0}, \beta)\ast \tdhat{d(\beta)}) &\cong \bigoplus_{\binom{d(\beta)}{l-i}} \redh_{l-d(\beta)-2}(\hat{0}, \beta) \\
&\cong \Z^{\binom{d(\beta)}{l-i}}\otimes \redh_{l-d(\beta)-2}(\hat{0}, \beta)\\
&\cong \redh^{l-i}(W_{\beta})\otimes \redh_{l-d(\beta)-2}(\hat{0}, \beta),
\end{align*}
which is the required decomposition.
\epr
Now we verify that the decomposition in Equation \ref{heq01} implies the result by De Concini and Procesi (Theorem \ref{thm0}).

\bc{cor0}
The Poincar\'e polynomial of the complement of a toric arrangement $\A$ is 
\[ \sum_{j=0}^l N_j (t+1)^{l-j} t^j.\]
\ec

\bpr 
Observe that the interval $[\hat{0}, \beta]$ is isomorphic to the intersection lattice of the local arrangement $\A[\beta]$. 
It follows from Example \ref{gmex02} that $\redh_{l-d(\beta)-2}(\hat{0}, \beta)\cong H^{l-d(\beta)}(M(\A[\beta]))$; hence its rank is a summand in $N_i$.
Finally, the other factor is the rank of $\redh^{l-i}(W_{\beta})$.
\epr

\be{linkex}
Consider the arrangement in Example \ref{ex3}. 
Here, $\widehat{B_2}\simeq S^1\vee S^2\vee S^2$. 
For each rank $1$ element there is a copy of $S^1\vee S^2$. 
For each rank $2$ element there is a copy of $S^1\vee S^1$. 
Hence, 
\[ \widehat{\ua} \simeq \bigvee_{10} S^1\vee \bigvee_{5} S^2.\]
The reduced cohomology of the complement is $\Z^{5}$ in degree $1$, $\Z^{10}$ in degree $2$ and zero everywhere else.
\ee

We now give an alternate description of the Goresky-MacPherson formula in terms of compactly supported cohomology. 
In a recent paper \cite{petersen17} Petersen constructed a spectral sequence for stratified spaces which computes the compactly supported cohomology of an open stratum.
The ingredients needed for this spectral sequence are the compactly supported cohomology groups of closed strata and the reduced cohomology groups of the poset of strata. 

Let $\A$ be a toric arrangement and for $\beta\in L_{\A}$ define $$W^{\circ}_{\beta} := W_{\beta}\setminus \bigcup_{\alpha>\beta}W_{\alpha}.$$
Consequently, we have a stratification $(\C^*)^l = \coprod_{\beta\in L_{\A}} W^{\circ}_{\beta}$. 
Note that $W^{\circ}_{\hat{0}} = M_{\A}$ and $\ol{W^{\circ}_{\beta}} = W_{\beta}$ for $\beta\in L_{\A}$. 
Since each stratum is an affine algebraic variety defined over $\C$ we have the following result by Petersen \cite[Theorem 3.3]{petersen17}. 
For notational simplicity denote by $L_p$ the set of codimension-$p$ components of $\A$.
\begin{thm}\label{petersenthm}
In the situation described above there is a spectral sequence of mixed Hodge structures
\[E_1^{p,q} = \bigoplus_{\beta\in L_p} \, \bigoplus_{i+j = p+q-2} \redh^i(\hat{0}, \beta)\otimes H^j_c(W_{\beta}) \Longrightarrow H^{p+q}_c(M_{\A}).\]
\end{thm}

Since in our case all the intervals in $L_{\A}$ are geometric lattices \cite[Example 3.10]{petersen17} the first page of the spectral sequence simplifies as follows:
\[ E_1^{p, q} = \bigoplus_{\beta\in L_p} \redh^{p-2}(\hat{0}, \beta)\otimes H^q_c(W_{\beta})\Longrightarrow H^{p+q}_c(M_{\A}).\]
These groups are nonzero only when $0\leq p\leq l$ and $l-p\leq q\leq 2l-2p$.
Moreover, all these columns have different weights so there can be no differentials in the spectral sequence.
Recall that for all $d\geq 0$ the group $H^k((\C^*)^d)$ is pure of weight $2k$ whereas the group $H^{2d-k}_c((\C^*)^d)$ is pure of weight $2d-2k$.
Consequently we have the following decomposition of the compactly supported cohomology of the complement. 

\bt{thmcpt}
Let $\A$ be a toric arrangement in $(\C^*)^l$ and $\ma$ be its complement. 
Then we have following analog of the Goresky-MacPherson formula using compactly supported cohomology:
\begin{align}
H^{l+j}_c(\ma) \cong \bigoplus_{p=0}^{l-j}\left(\bigoplus_{\beta\in L_p} \redh^{p-2}(\hat{0}, \beta)\otimes H_c^{l+j-p}(W_{\beta})\right).
\end{align}
\et 
Note that the above decomposition is Poincar\'e dual to the earlier one in Equation \ref{heq01}. 
\be{cptex1}
We consider the toric arrangement in Example \ref{ex3}. 
Here, an order complex of intervals ending in rank $2$ elements consists of $3$ discrete points hence its reduce cohomology is $\Z^2$ in dimension $0$. 
An easy calculation shows that for $i= 2, 3, 4$  $H^i_c(\ma) = \Z^{10}, \Z^5, \Z$ respectively and $0$ for all other values of $i$. 
\ee
\be{cptex2}
Now consider the toric arrangement in Example \ref{ex4}. 
For a space $X$ denote by $\mathcal{F}(X ,n)$ the configuration space of $n$ ordered points in $X$.
Using above result and Poincar\'e duality we get
\begin{align*}
H_i(\mathcal{F}(\C^* ,n)) &\cong \bigoplus_{p=0}^i\left(\bigoplus_{\beta\in(\Pi_n)_p} \redh^{p-2}(\hat{0}, \beta) \otimes H_{i-p}(W_\beta)\right).
\end{align*}
\ee

\section{Concluding remarks}
In the final section we consider collections of arbitrary toric subspaces, i.e., they could have codimension greater than $1$.
The definition of intersection poset is the exact same as in the codimension $1$ case and we denote this poset again by $L_{\A}$. 
In general, this poset need not be graded and the intervals in it need not be geometric lattices.
The complement of a toric subspace arrangement is also defined analogously and is denoted by $\ma$.
To the best of our knowledge these kinds of arrangements have not been considered before.
We first state the Goresky-MacPherson formula in this general context and then list some future research directions. 

\begin{thm}\label{mainthm}
Suppose that $\A$ is a finite collection of toric subspaces not necessarily of codimension-$1$ and $L_{\A}$ is the associated poset of intersections. 
Then for every $k\geq 0$ we have
\begin{align}
\tld{H}^k(\ma) \cong 
\bigoplus_{\beta\in L_{\A}}
\left[\bigoplus_{j=0}^{d(\beta)}\left(
\bigoplus_{\binom{d(\beta)}{j}}\redh_{2l-2-j-k-d(\beta)}(\hat{0}, \beta)
\right)
\right].\label{heq03}
\end{align}
\end{thm}

\bpr 
The modified link $\ua$ of $\A$ is the union of all coordinate axes and the toric subspaces in $\A$.  
Note that the expression for the homotopy type of $\widehat{\ua}$, its one point compactification, remains the same as in Theorem \ref{thm03}.
Application of Alexander duality calculations similar to the ones in the proof of Theorem \ref{toricgm} give the above decomposition. 
\epr

We also derive an expression for the Euler characteristic which exhibits the role of the M\"obius function.

\bt{eulerc}
The Euler characteristic of $M_{\A}$ is determined by $0$-dimensional intersections as given by the expression
\begin{align}\label{euexp}
    \chi(M_{\A}) = \sum_{d(\beta) = 0} \mu(\hat{0}, \beta).
\end{align}
\et

\bpr 
It should be clear that computations using Betti numbers are tedious so employ a different strategy. 
Since $\ma$ is a non-compact, $2l$-manifold $\chi_c(\ma) = \chi(\ma)$.
Recall that the compactly supported Euler characteristic is additive over stratifications.
Using the stratification in Theorem \ref{thmcpt}
we get $\chi_c((\C^{\ast})^l) = \sum_{\beta\in L}\chi_c(W_{\beta}^{\circ})$.
By the M\"obius inversion formula for the poset $L$, we have
\[\chi_c(\ma) = \sum_{\beta\in L}\mu(\hat{0}, \beta) \cdot\chi_c(W_{\beta}). \]
Now use the fact that $\chi_c((\C^{\ast})^d) = 0$ whenever $d\geq 1$.
\epr

If the intersection poset is (sequentially) Cohen-Macaulay then the cohomology of the complement is torsion free. 
For example, one could consider arrangements $\A$ such that all the members of $\A$ are of fixed codimension. 
One can also consider $k$-equal arrangements in the toric context, i.e., collections of subspaces with at least $k$ identical coordinates. 
For this class of arrangements the intersection poset is the $k$-equal partition lattice; which is known to be sequentially Cohen-MaCaulay \cite[Section 3.2.4]{wachs07}. 

Analogous to the linear case one can also consider ``coordinate subspace arrangements'' defined by simplicial complexes.
For example, if $\Delta$ is a simplicial complex then consider $\A = \{\{z_i = 1\mid i\in \sigma\}\mid \sigma\in \Delta\}$. 
The intersection poset of this arrangement is anti-isomorphic to the face poset of $\Delta$.
A Choice of an appropriate $\Delta$ will introduce torsion in $H^i(\ma, \Z)$.

It should be interesting to see if there is a combinatorial realization of cohomology classes. 
One could also try to figure out a combinatorial description of the multiplication in the associated graded? (See \cite[Section 11]{vas01} for the linear case).
A more challenging problem could be to find a presentation for the cohomology algebra ${H}^{\ast}(\ma, \Z)$. 

Another interesting property of linear subspace arrangements is that their complements are formal. 
This is also known to be true for arrangement of toric hypersurfaces. 
What happens for a general toric subspace arrangement?

In the linear subspace case it was proved that the stable homotopy type of the complement is combinatorially determined, see \cite[Theorem 3.4]{zz93} or \cite[Corollary 4]{vas01}. 
So given a toric arrangement $\A$ one can try and determine a number $N$ such that the hmotopy type of iterated suspension $\Sigma^N\ma$ stabilizes.
Moreover, is it true that the (single) suspension of the complement of hyper-toric arrangement (i.e., all members of $
\A$ are of codimension $1$) is a wedge of spheres?
This is the case for complex hyperplane arrangements as proved by Schaper in \cite{schaper97}.

Consider the arrangement $\A = \{z_i = 1, z_jz_k^{\pm 1} =1 \mid 1\leq i,j,k\leq l, j\neq k\}$ in $(\C^{\ast})^l$.
Applying the biholomorphic map $z\mapsto \frac{1+z}{1-z}$ coordinate wise to the complement of $\A$ we get that $\ma\cong\{ w\in \C^l\mid w_i\neq \pm w_j, w_i\neq \pm 1\}$.
This is a complement of a hyperplane arrangement. 
There are toric arrangements whose complements have the homotopy type of a complement of a hyperplane arrangement. 
For example, $\mathcal{F}(\C, n)\stackrel{\simeq}{\to} \mathcal{F}(\C^{\ast}, n-1)$ via $z_i\mapsto z_i-z_n$.
Hence, one can ask for a qualitative characterization of those toric (subspace) arrangements whose complement has the homotopy type of the complement of linear subspace arrangement.

\bibliographystyle{abbrv} 
\bibliography{deshpande_gm_ref} 

\end{document}